\documentclass[12pt]{article}
\usepackage[cp1251]{inputenc}
\usepackage[russian]{babel}
\usepackage{amssymb}
\usepackage{amsmath}

\begin{document}

\begin{center}
\textbf{\large On the Fourier transform of the characteristic
functions of domains with $C^1$ -smooth boundary}
\end{center}

\begin{center}
Vladimir Lebedev
\end{center}

\begin{quotation}
\small{\textsc{Abstract.} We consider domains $D\subseteq\mathbb
R^n$ with $C^1$ -smooth boundary and study the following question:
when the Fourier transform $\widehat{1_D}$ of the characteristic
function $1_D$ belongs to $L^p(\mathbb R^n)$?

  References: 16 items.


  }
\end{quotation}

\begin{center}
\textbf{Introduction}
\end{center}

   Let $D$ be a bounded domain (an open and connected set)
in $\mathbb R^n, ~n\geq 2$. Consider its characteristic function
$1_D$, i.e. the function that takes value $1_D(t)=1$ for $t\in D$
and value $1_D(t)=0$ for $t\notin D$. Consider the Fourier
transform $\widehat{1_D}$ of this function. In the present work we
study the following question: for which domains we have
$\widehat{1_D}\in L^p(\mathbb R^n)$? Only the case when $1<p<2$ is
interesting.

  It will be convenient for us to deal with the spaces
$A_p(\mathbb R^n), ~1\leq p\leq\infty,$ of tempered distributions
$f$ on $\mathbb R^n$ such that the Fourier transform $\widehat{f}$
belongs to $L^p(\mathbb R^n)$. The norm on $A_p(\mathbb R^n)$ is
defined in the natural way:
$$
\|f\|_{A_p(\mathbb R^n)}=\|\widehat{f}\|_{L^p(\mathbb R^n)}.
$$
Recall that (see, e.g. [1, Ch. V, \S~1]) for $1\leq p\leq 2$ the
Fourier transform (as well as its inverse) is a bounded operator
from $L^p(\mathbb R^n)$ to $L^q(\mathbb R^n), ~1/p+1/q=1$, so each
distribution in $A_p(\mathbb R^n), ~1\leq p\leq 2,$ is a function
in $L^q(\mathbb R^n)$.

  Direct calculation shows that if $D$ is a cube in $\mathbb
R^n$, then $1_D\in A_p(\mathbb R^n)$ for all $p>1$. The same is
true in the case when $D$ is a polytope (i.e. a finite union of
simplices). On the other hand, using the well known asymptotics
for Bessel functions, one can verify that if $D\subseteq\mathbb
R^n$ is a ball, then $1_D\in A_p(\mathbb R^n)$ for $p>2n/(n+1)$
and $1_D\notin A_p(\mathbb R^n)$ for $p\leq 2n/(n+1)$. The same
result holds in the general case of bounded domains with twice
smooth boundary. (This follows from Theorems 1, 2 of the present
work, see Corollary 2.) Thus, for (bounded) domains with $C^2$
-smooth boundary $2n/(n+1)$ is the critical power of integrability
for the Fourier transform of the characteristic function.

    In the present work we shall obtain a series of results on the
behavior of the Fourier transform  of the characteristic functions
of bounded domains with $C^1$ -smooth boundary. Generally
speaking, this case is essentially different from the twice smooth
case as is shown (see \S~3) by an example of a domain
$D\subseteq\mathbb R^2$ whose boundary is $C^1$ -smooth and at the
same time $1_D\in A_p(\mathbb R^2)$ for all $p>1$. (The critical
value for planar domains with twice smooth boundary is $4/3$.)

  We note that various questions on the rate of decrease
at infinity for the Fourier transform of characteristic functions
of domains and closely related questions on the behavior of the
Fourier transform of (smooth) measures supported on surfaces were
investigated by many authors and represent a classical topic in
harmonic analysis, see Stein's survey [2], where one can find an
ample bibliography, and his book [3] (Ch. VIII). The basic tools
to obtain asymptotic estimates in these investigations are the
stationary phase method and the van der Corput lemma. The use of
this methods requires considerable smoothness of the boundary of a
domain. The smoothness should be at least two even in the planar
case. The crucial role in this approach is played by the curvature
of the surface (of the boundary of a domain). Our approach does
not use any arguments related to curvature and allows to consider
domains with $C^1$ -smooth boundary.

   We denote by $\partial D$ the boundary of a domain
$D\subseteq\mathbb R^n$. Saying that the boundary of $D$ is $C^1$
-smooth or $C^2$ -smooth we mean that in an appropriate
neighborhood of each of its points the boundary $\partial D$ is a
graph of a certain (real) function of class $C^1$ or $C^2$
respectively (that is of a function whose all partial derivatives
of the first or the second order respectively are continuous).

  For each domain $D\subseteq\mathbb R^n$ with $C^1$ -smooth
boundary let $\nu_D (x)$ be the outer unit normal vector to
$\partial D$ at a point $x\in\partial D$. The corresponding map
$\nu_D:\partial D \rightarrow S^{n-1}$ of the boundary of $D$ into
the unit sphere $S^{n-1}$ centered at the origin is called normal
map. By $\omega (\nu_D, \delta)$ we denote the modulus of
continuity of $\nu_D$:
$$
\omega (\nu_D, \delta)=\sup_{x, y ~\in \partial D; ~|x-y|\leq\delta}
|\nu_D(x)-\nu_D(y)|, \qquad \delta\geq 0,
$$
where $|u|$ is the length of a vector $u\in\mathbb R^n$. Let then
$\omega (\delta)$ be an arbitrary nondecreasing continuous
function on $[0, \infty)$, $\omega (0)=0$. In the case when
$\omega (\nu_D, \delta)=O(\omega (\delta)), ~\delta\rightarrow
+0,$ we say that the boundary $\partial D$ is $C^{1, \omega}$
-smooth. For bounded domains this condition is equivalent to the
condition that in an appropriate neighborhood of each of its
points the boundary of $D$ is a graph of a certain function of
class $C^{1, \omega}$. In other words, for each point
$x\in\partial D$ one can find a neighborhood $B$, containing $x$,
and a domain $V\subseteq\mathbb R^{n-1}$ such that $B\cap\partial
D$ is a graph of some (real) function $\varphi\in C^{1,
\omega}(V)$ i.e. of a function with $\omega(V,\nabla\varphi,
\delta)=O(\omega(\delta)), ~\delta\rightarrow +0$, where
$$
\omega(V,\nabla\varphi, \delta)=\sup_{x, y\in V; |x-y|\leq
\delta}|\nabla\varphi(x)-\nabla\varphi(y)|, \qquad \delta\geq 0,
$$
is the modulus of continuity of the gradient $\nabla\varphi$ of
$\varphi$.

  If the boundary $\partial D$ of a domain $D$ is $C^1, ~C^2$, or
$C^{1, \omega}$ -smooth, we write $\partial D\in C^1, ~\partial
D\in C^2$, or $\partial D\in C^{1, \omega}$ respectively.

  If $\omega (\delta)=\delta^\alpha, ~0<\alpha\leq 1$, then we
just write $C^{1, \alpha}$ instead of $C^{1, \delta^\alpha}$.

  In \S ~1 we give a simple proof of the inclusion
$1_D\in A_p(\mathbb R^n)$, which holds for $p>2n/(n+1)$ for all
bounded domains $D\subseteq\mathbb R^n$ with $C^1$ -smooth
boundary (Theorem 1). For convex domains (without smoothness
assumptions on the boundary) such assertion was earlier obtained
by Herz [4].

  In \S ~2 we obtain the main result of the present work.
Namely, we show (Theorem 2) that if $\partial D\in C^{1, \omega}$
and
$$
\int_0^1 \frac{\delta^{n(p-1)-1}}{\omega(\delta)^{n-p}} d\delta
=\infty,
$$
then $1_D\notin A_p(\mathbb R^n)$. In particular (Corollary 1), if
$\partial D\in C^{1, \alpha}$, then $1_D\notin A_p(\mathbb R^n)$
for
$$
p\leq 1+\frac{(n-1)\alpha}{n+\alpha}.
$$
Putting $\alpha=1$ here and taking into account the preceding
result, we obtain the indicated in the beginning of the
introduction assertion on the critical value for domains with
twice smooth boundary (Corollary 2).

  In \S ~3 we consider planar domains. According to the result
indicated above if for a domain $D\subseteq\mathbb R^2$ we have
$\partial D\in C^{1, \omega}$ and
$$
\int_0^1 \frac{\delta^{2p-3}}{\omega(\delta)^{2-p}} d\delta
=\infty,
$$
then $1_D\notin A_p(\mathbb R^2)$. In particular, if $\partial
D\in C^{1, \alpha}$, then $1_D\notin A_p(\mathbb R^2)$ for $p\leq
1+\alpha/(2+\alpha)$. We show (Theorem 3) that this result is
sharp, namely, for each class $C^{1, \omega}$ (under certain
simple condition imposed on $\omega$) there exists a bounded
domain $D\subseteq\mathbb R^2$ such that $\partial D\in C^{1,
\omega}$ and for all $p>1$ satisfying
$$
\int_0^1 \frac{\delta^{2p-3}}{\omega(\delta)^{2-p}} d\delta
<\infty
$$
we have $1_D\in A_p(\mathbb R^2)$. In particular (Corollary 3), if
$0<\alpha<1$, then there exists a planar domain $D$ with $C^{1,
\alpha}$ -smooth boundary such that $1_D\in A_p(\mathbb R^2)$ for
all $p>1+\alpha/(2+\alpha)$. It also follows that (Corollary 4)
there exists a planar domain $D$ with $C^1$ -smooth boundary such
that $1_D\in A_p(\mathbb R^2)$ for all $p>1$ (it suffices to take
$\omega$ decreasing to zero slower than any power, i.e. so that
$\lim_{\delta\rightarrow +0}\omega
(\delta)/\delta^\varepsilon=\infty$ for all $\varepsilon>0$).

   The results of the present work are essentially based on the
results obtained by the author in [5], [6] where for real $C^1$
-smooth functions $\varphi$ we studied the growth of the $A_p$
norms of exponential functions $e^{i\lambda\varphi}$. The question
on the growth of the norms of these functions naturally arises in
relation with the known Beurling--Helson theorem (see the history
of the question in [5]). Simple arguments (Lemma 1 of the present
work) allow to reduce the study of the characteristic functions to
the study of the behavior of the exponential functions.

   We denote by $|E|$ the Lebesgue measure of a (measurable) set
$E\subseteq\mathbb R^n$ and by $|E|_{S^{n-1}}$ the spherical
measure of a set $E\subseteq S^{n-1}$. We use $(x, y)$ to denote
the inner product of vectors $x$ and $y$ in $\mathbb R^n$. If
$E\subseteq\mathbb R^n$ and $t\in\mathbb R^n$, then we put
$E+t=\{x+t: x\in E\}$. Various positive constants are denoted by
$c, c_p, c_{p, n}$.

\quad

  The results of this work were partially presented at
11-th and 14-th Summer St. Petersburg Meetings in Mathematical
Analysis [7], [8] and completely at the International Conference
``Harmonic Analysis and Approxi-mations, III'', Tsahkadzor,
(Armenia) [9].

\quad

  The author is grateful to E. A. Gorin who read the introduction
and made useful remarks.

\quad

\begin{center}
\textbf{\S ~1. General case of domains with $C^1$ -smooth
boundary}
\end{center}

\textbf{Theorem 1.} \emph{Let $D$ be a bounded domain in $\mathbb
R^n, ~n\geq 2,$ with $\partial D\in C^1$. Then $1_D\in A_p(\mathbb
R^n)$ for all $p>2n/(n+1)$.}

\quad

\emph{Proof.} For $s>0$ consider the Sobolev spaces $W_2^s(\mathbb
R^n)$ of functions $f\in L^2 (\mathbb R^n)$ satisfying
$$
\|f\|_{W_2^s(\mathbb R^n)}=\Big(\int_{\mathbb R^n}
(|\xi|^{2s}+1)|\widehat{f}(\xi)|^2 d\xi \Big)^{1/2} <\infty.
$$

   It is easy to verify that
$$
W_2^s(\mathbb R^n)\subseteq A_p(\mathbb R^n)
$$
for $2n/(n+2s)<p<2$. Indeed, put $p^\ast=2/p,
~1/p^\ast+1/q^\ast=1$. We have $spq^\ast>n$. Using the H\"older
inequality, we obtain
\begin{multline*}
\|f\|_{A_p(\mathbb R^n)}^p=\int_{\mathbb R^n} |\widehat{f}(\xi)|^p
d\xi= \int_{\mathbb R^n}
\Big(|\widehat{f}(\xi)|\,(|\xi|^s+1)\Big)^p
\frac{1}{(|\xi|^s+1)^p} d\xi \\ \leq \bigg(\int_{\mathbb R^n}
\Big(|\widehat{f}(\xi)|\,(|\xi|^s+1)\Big)^{pp^\ast}
d\xi\bigg)^{1/p^\ast}\bigg(\int_{\mathbb R^n}
\frac{1}{(|\xi|^s+1)^{pq^\ast}} d\xi\bigg)^{1/q^\ast}\leq c_{p, s}
\|f\|_{W_2^s (\mathbb R^n)}^p.
\end{multline*}

  To prove the theorem it remains to take into account that for
each bounded domain $D\subseteq\mathbb R^n$ with $C^1$ -smooth
boundary we have $1_D\in W_2^s(\mathbb R^n)$ for all $s<1/2$. This
is a trivial consequence of the theorem on (pointwise) multipliers
of the Sobolev spaces [10, \S~5].

  We shall give an independent short and simple proof of the inclusion
$1_D\in W_2^s, ~s<1/2$. It is well known that for $0<s<1$ the norm
$\|\cdot\|_{W_2^s(\mathbb R^n)}$ and the norm
$$
\|f\|=\|f\|_{L^2(\mathbb R^n)}+\bigg(\int_{\mathbb R^n}
\frac{1}{|t|^{n+2s}}
\bigg(\int_{\mathbb R^n}|f(x+t)-f(x)|^2 dx\bigg)dt\bigg)^{1/2}
\eqno(1)
$$
are equivalent (see, e.g., [11, Ch. V, \S ~3.5 ]). Note now that
for each $t\in\mathbb R^n$ the symmetric difference
$$
((D-t)\setminus D)\cup (D\setminus (D-t))
$$
of the sets $D-t$ and $D$ is contained in the (closed) $|t|$
-neighborhood of the boundary $\partial D$ of $D$, so its
(Lebesgue) measure is at most $c |t|$. It is also clear that the
measure of this symmetric difference is at most $2|D|$. Thus,
$$
\int_{\mathbb R^n}|1_D(x+t)-1_D(x)|^2 dx \leq \min (c |t|, 2|D|),
\qquad t\in\mathbb R^n,
$$
and it remains to use the equivalence of the norm
$\|\cdot\|_{W_2^s(\mathbb R^n)}$ and the norm defined by (1). The
theorem is proved. \footnote{Note that [12, Corollary 2.2] if
$E\subseteq \mathbb R^n, ~n\geq 1,$ is a set of positive measure,
then $1_E\notin W_2^{1/2}(\mathbb R^n)$.}

\quad

  \emph{Remark} 1. The method used in the proof of Theorem 1
can be applied to arbitrary sets (not only domains). Recall [13]
that the upper Minkowski dimension $\overline{\dim}_M F$ of a
bounded set $F\subseteq\mathbb R^n$ is defined by
$$
\overline{\dim}_M F=\inf\{0\leq\gamma\leq n : |(F)_\delta|=
O(\delta^{n-\gamma}), ~\delta\rightarrow+0\},
$$
where $(F)_\delta$ is the $\delta$ -neighborhood of $F$. Let
$E\subseteq\mathbb R^n, ~n\geq 1,$ be a bounded set of positive
measure. Let $a$ be the upper Minkowski dimension of its boundary
$\partial E$. Suppose that $a<n$. Repeating with obvious
modifications the arguments used above, we see that then $1_E\in
W_2^s(\mathbb R^n)$ for all $s<(n-a)/2$. Hence in turn we obtain
that $1_E\in A_p(\mathbb R^n)$ for all $p>2n/(2n-a)$. Note that
for $s<(n-a)/2$ the usage of the norm (1) is justified since
$a\geq n-1$. Indeed, let us verify that if a set
$E\subseteq\mathbb R^n$ is bounded and has positive measure, then
$\overline{\dim}_M \partial E\geq n-1$. Assuming that
$E\setminus\partial E\neq\varnothing$ (otherwise there is nothing
to prove) fix a point $x_0\in E\setminus\partial E$. There exists
an open ball $B$ centered at $x_0$ that does not contain points of
the boundary $\partial E$ and moreover lies at positive distance
from $\partial E$. Denote by $S$ the boundary sphere of the ball
$B$. Define a map $\theta : \mathbb R^n\setminus B\rightarrow S$
as follows. Take a point $x\in \mathbb R^n\setminus B$ and
consider the ray that passes trough $x$ and has its origin at
$x_0$. Denote by $\theta(x)$ the point of intersection of this ray
with the sphere $S$. Clearly the map $\theta$ is Lipschitz
(moreover it is non-expanding, i.e. $|\theta(x_1)-\theta(x_2)|\leq
|x_1-x_2|$ for all $x_1, x_2 \in \mathbb R^n\setminus B$). It is
easy to see that the image of the boundary of the set $E$ under
the map $\theta$ is the whole sphere $S$. At the same time it is
known [13, Ch. 7] that Lipschitz maps do not increase the
dimension of a set. Thus,
$$
n-1=\overline{\dim}_M S=\overline{\dim}_M \theta(\partial E)\leq
\overline{\dim}_M \partial E.
$$

\quad

\begin{center}
\textbf{\S ~2. Domains with $C^{1, \omega}$ -smooth boundary}
\end{center}

\textbf{Theorem 2.} \emph{Let $D$ be a bounded domain in $\mathbb
R^n, ~n\geq 2,$ with $\partial D\in C^{1, \omega}$. If
$$
\int_0^1\frac{\delta^{n(p-1)-1}}{(\omega(\delta))^{n-p}}~d\delta=\infty,
$$
then $1_D\notin A_p(\mathbb R^n)$.}

\quad

  From Theorem 2 we immediately obtain the following corollary.

\quad

\textbf{Corollary 1.} \emph{Let $0<\alpha\leq 1$. Let $D$ be a
bounded domain in $\mathbb R^n, ~n\geq 2,$ with $\partial D\in
C^{1, \alpha}$. If
$$
p\leq 1+\frac{(n-1)\alpha}{n+\alpha},
$$
then $1_D\notin A_p(\mathbb R^n)$.}

\quad

  Note particularly the case of domains with twice smooth
boundary and even more general $C^{1, 1}$ case. Namely, using
Corollary 1 and Theorem 1, we obtain the following corollary.

\quad

  \textbf{Corollary 2.}  \emph{Let $D$ be a bounded domain in
$\mathbb R^n, ~n\geq 2,$ with $\partial D\in C^{1, 1}$. Then
$1_D\in A_p(\mathbb R^n)$ for $p>2n/(n+1)$ and $1_D\notin
A_p(\mathbb R^n)$ for $p\leq 2n/(n+1)$. In particular, this holds
for bounded domains with twice smooth boundary.}

\quad

  Before proving the theorem we discuss certain preliminaries
and prove a certain lemma.

  Recall (see e.g. [11, Ch. IV, \S~3.1]) that a function
$m\in L^\infty(\mathbb R^n)$ is called an $L^p$ -Fourier
multiplier ($1\leq p\leq\infty$) if the operator $Q$ given by
$$
\widehat{Qf}=m\widehat{f}, \qquad f\in L^p\cap L^2(\mathbb R^n),
$$
is a bounded operator from $L^p(\mathbb R^n)$ to itself. The space
$M_p(\mathbb R^n)$ of all such multipliers endowed with the norm
$$
\|m\|_{M_p(\mathbb R^n)}=
\|Q\|_{L^p(\mathbb R^n)\rightarrow L^p(\mathbb R^n)}
$$
is a Banach algebra (with respect to the usual multiplication of
functions). It is well known that the characteristic function of
an arbitrary parallelepiped is a multiplier for all $p,
1<p<\infty,$ (see e.g. [11, Ch. IV, \S~4.1] and also [14, Ch. I,
\S~1.3]).

  Note that in the case when $1<p<2$ (only this case
is interesting to us here) each $L^p$ -Fourier multiplier is a
pointwise multiplier of $A_p$, that is if $m\in M_p(\mathbb R^n)$,
then for every function $f\in A_p(\mathbb R^n)$ we have $m f\in
A_p(\mathbb R^n)$ and
$$
\|m f\|_{A_p}\leq \|m\|_{M_p}\|f\|_{A_p}.
\eqno(2)
$$
This can be easily verified as follows. Note that estimate (2)
holds for every function $f\in A_p\cap L^2(\mathbb R^n)$ (this is
obvious since the Fourier transform and its inverse differ only in
sign of variable and normalization factor). It is clear that the
set $A_p\cap L^2(\mathbb R^n)$ is dense in $A_p(\mathbb R^n)$. Let
$f\in A_p(\mathbb R^n)$ be an arbitrary function. Let $f_k\in
A_p\cap L^2(\mathbb R^n), ~k=1, 2, \ldots,$ be a sequence that
converges to $f$ in $A_p(\mathbb R^n)$. We have
$$
\|mf_j-mf_k\|_{A_p}=\|m\cdot (f_j-f_k)\|_{A_p}\leq
\|m\|_{M_p}\|f_j-f_k\|_{A_p}\rightarrow 0.
$$
It is obvious that the spaces $A_p(\mathbb R^n)$ are Banach
spaces, so the sequence $m f_k, ~k=1, 2, \ldots,$ converges in
$A_p(\mathbb R^n)$ to some function $g\in A_p(\mathbb R^n)$. Using
the Hausdorff--Young inequality [1, Ch. V, \S~1], which in our
notation has the form $\|\cdot\|_{L^q}\leq\|\cdot\|_{A_p},
~1/p+1/q=1, ~1\leq p\leq 2,$ we see that the sequences $\{f_k\}$
and $\{m f_k\}$ converge in $L^q(\mathbb R^n)$ to $f$ and $g$
respectively. Thus, $m f=g$ and it remains to proceed to the limit
in the inequality $\|m f_k\|_{A_p}\leq \|m\|_{M_p}\|f_k\|_{A_p}$.

  Let $D_1$ be a domain in $\mathbb R^n$ and let
$D_2=l(D_1)$ be its image under a non-degenerate affine map
$l:\mathbb R^n\rightarrow\mathbb R^n$. It is easy to see that
$1_{D_1}\in A_p(\mathbb R^n)$ if and only if $1_{D_2}\in
A_p(\mathbb R^n)$. It suffices to observe that if $l(x)=Qx+b$,
then for each function $f\in L^1(\mathbb R^n)$ we have
$|\widehat{f\circ l}(u)|=|\det Q|^{-1}|\widehat{f}((Q^{-1})^*u)|$,
where $Q^{-1}$ is the matrix inverse of $Q$ and $(Q^{-1})^*$ is
the matrix adjoint to $Q^{-1}$.

  Let $E$ be an arbitrary set in $\mathbb R^m, ~m\geq 1$.
Following [6] we say that a function $f$ defined on $E$ belongs to
the space $A_p(\mathbb R^m, E)$ if there exists a function $F\in
A_p(\mathbb R^m)$ such that its restriction $F_{| E}$ to the set
$E$ coincides with $f$. We define the norm on $A_p(\mathbb R^m,
E)$ by
$$
\|f\|_{A_p(\mathbb R^m, E)}=\inf_{F_{| E}=f} \|F\|_{A_p(\mathbb R^m)}.
$$

  Note that if $I$ is a parallelepiped in $\mathbb R^m$ and $f$
is a function on $I$, then putting
$$
\|f\|_{A_p(\mathbb R^m, I)}^\circ =\|F\|_{A_p(\mathbb R^m)},
\eqno(3)
$$
where $F$ is the function $f$ extended by zero to the compliment
$\mathbb R^m\setminus I$ (that is $F=f$ on $I$ and $F=0$ on
$\mathbb R^m\setminus I$), we obtain the norm
$\|\cdot\|_{A_p(\mathbb R^m, I)}^\circ $ equivalent to the norm
$\|\cdot\|_{A_p(\mathbb R^m, I)}$ for $1<p<2$. This follows since
for $1<p<\infty$ the characteristic function of a parallelepiped
is an
 $L^p$ -Fourier multiplier.

  The following result, obtained by the author in
[6, Theorem 1$''$], is the base of the proof of Theorem 2. Let
$1\leq p <2$. Let $V$ be a domain in $\mathbb R^m$ and let
$\varphi\in C^{1,\omega}(V)$ be a real function. Suppose that the
gradient $\nabla\varphi$ of $\varphi$ is non-degenerate on $V$,
i.e. the set $\nabla\varphi (V)$ is of positive measure. Then for
all those $\lambda\in\mathbb R, ~|\lambda|\geq 1,$ for which
$e^{i\lambda\varphi}\in A_p(\mathbb R^m, V)$, we have
$$
\|e^{i\lambda\varphi}\|_{A_p(\mathbb R^m, V)} \geq
c \bigg(|\lambda|^{1/p}\chi^{-1}\bigg(\frac{1}{|\lambda|}\bigg)\bigg)^m,
\eqno(4)
$$
where $\chi^{-1}$ is the function inverse to
$\chi(\delta)=\delta\omega(\delta)$ and $c=c(p, \varphi)>0$ is
independent of $\lambda$.

   Simple Lemma 1 below (which we shall also use
in \S ~3) allows to reduce the question on inclusion $1_D\in A_p$
to the question on behavior of exponential functions
$e^{i\lambda\varphi}$ in $A_p$.

  For a vector $x=(x_1, x_2, \ldots, x_{m})\in\mathbb
R^{m}$ and a number $a\in \mathbb R$ let $(x, a)$ denote the
vector $(x_1, x_2, \ldots, x_{m}, a)\in\mathbb R^{m+1}$.

  Let $I$ be an open parallelepiped in $\mathbb R^m$ with
edges parallel to coordinate axes. Let $\varphi$ be a continuous
bounded function on $I$ such that $\varphi (t)>0$ for all $t\in
I$. Consider the following domain $G$ in $\mathbb R^{m+1}$
$$
G=\{(t,y)\in \mathbb R^m\times \mathbb R : t\in I, ~0<y<\varphi
(t)\}.
$$
Each domain of this form is called a special domain generated by
the pair $(I, \varphi)$.

\quad

\textbf{Lemma 1.} \emph{Let $G\subseteq\mathbb R^{m+1}$ be a
special domain generated by a pair $(I, \varphi)$. Let $1<p<2$.
The inclusion $1_G\in A_p(\mathbb R^{m+1})$ holds if and only if
$e^{i\lambda\varphi}\in A_p(\mathbb R^m, I)$ for almost all
$\lambda\in\mathbb R$ and}
$$
\int_{\mathbb R} \frac{1}{|\lambda|^p}
\|e^{i\lambda\varphi}-1\|_{A_p(\mathbb R^m, I)}^p ~d\lambda <\infty.
$$

\quad

\emph{Proof.}  For $\lambda\in\mathbb R\setminus\{0\}$ define the
function $F_\lambda$ on $\mathbb R^m$ by
$$
F_\lambda(t)=\left\{ \begin{array}{ll}
\displaystyle\frac{1}{-i\lambda}(e^{-i\lambda\varphi
(t)}-1), &\textrm{if} ~~t\in I,\\
0, &\textrm{if} ~~t\in \mathbb R^m\setminus I. \end{array}\right.
$$

   Note that
$$
\widehat{1_G}(u, \lambda)=\widehat{F_\lambda}(u),
\qquad (u, \lambda)\in \mathbb
R^m\times\mathbb R, \quad \lambda\neq 0.
$$
Indeed, direct calculation yields
$$
\widehat{1_G}(u, \lambda)=\underset{t\in I ~0<y<\varphi (t)}
{\int\int}e^{-i(u, t)} e^{-i\lambda y}dtdy=\int_I \bigg (\int_0^{\varphi
(t)}e^{-i\lambda y}dy\bigg )e^{-i(u, t)}dt
$$
$$
=\int_I \frac{1}{-i\lambda}(e^{-i\lambda\varphi (t)}-1) e^{-i(u,
t)}dt=\widehat{F_\lambda}(u).
$$

  Thus, $1_G\in A_p(\mathbb R^{m+1})$ if and only if
$$
\int_{\mathbb R} \|F_\lambda\|_{A_p(\mathbb R^m)}^p d\lambda<\infty.
$$

  It remains only to take into account that
$$
\|F_\lambda\|_{A_p(\mathbb R^m)}=\frac{1}{|\lambda|}
\|e^{i\lambda\varphi}-1\|_{A_p(\mathbb R^m, I)}^\circ,
$$
where $\|\cdot\|_{A_p(\mathbb R^m, I)}^\circ$ is the equivalent
norm on $A_p(\mathbb R^m, I)$ defined above (see (3)). The lemma
is proved.

\quad

  \emph{Proof of Theorem} 2. Assume that contrary to the
assertion of the theorem we have $1_D\in A_p(\mathbb R^n)$.

   Each point $x$ of the boundary $\partial D$ can be surrounded
by an open parallelepiped $\Pi_x\ni x$ so small that after an
appropriate rotation and translation the intersection $D\cap\Pi_x$
becomes a special domain. Extract a finite subcovering from the
covering $\{\Pi_x, ~x\in \partial D\}$ of $\partial D$. Note that
$\nu_D(\partial D)=S^{n-1}$, so at least for one of the
parallelepipeds $\Pi_x$, which we denote by $\Pi$, we have
$$
|\nu_D (\partial D\cap \Pi)|_{S^{n-1}}>0. \eqno(5)
$$
Consider the domain $G=D\cap \Pi$. Since the characteristic
function of a parallelepiped is an $L^p$ -multiplier, we have
$1_G=1_\Pi \cdot 1_D\in A_p(\mathbb R^n)$. Replacing, if needed,
the domain $G$ by its copy obtained by rotation and translation,
we can assume that $G$ is a special domain. This domain is
generated by a pair $(I, \varphi)$, where $I$ is a certain
parallelepiped in $\mathbb R^m, ~m+1=n,$ with edges parallel to
coordinate axes and $\varphi$ is a certain function in $C^{1,
\omega}(I)$.

   It is easy to see that condition (5) implies that the gradient
of $\varphi$ is non-degenerate on $I$, that is, we have
$|\nabla\varphi(I)|>0$. Indeed (recall that if $x=(x_1, x_2,
\ldots, x_m)\in\mathbb R^m$ and $a\in \mathbb R$, then $(x, a)$
denotes the vector $(x_1, x_2, \ldots, x_m, a)\in\mathbb R^{m+1}$)
consider the map
$$
\beta(t)=(t, \varphi (t)), \qquad t\in I,
$$
($\beta$ maps $I$ onto the graph of $\varphi$). The normal map
$\nu_D$ and the gradient $\nabla\varphi$ of $\varphi$ are related
by
$$
\nu_D\circ\beta(t)=
\frac{1}{\sqrt{|\nabla\varphi(t)|^2+1}}(-\nabla\varphi(t), 1),
\qquad t\in I.
$$
Thus, putting
$$
\gamma (\xi)=\frac{1}{\sqrt{|\xi|^2+1}}(-\xi, 1), \qquad \xi\in
\mathbb R^m,
$$
we have $\nu_D\circ\beta=\gamma\circ\nabla\varphi$. So for the set
$W=\bigtriangledown\varphi(I)$ we obtain
$$
\gamma(W)=\gamma(\nabla\varphi(I))=\nu_D\circ\beta(I)=\nu_D(\partial
D\cap \Pi),
$$
and relation (5) implies $|\gamma(W)|_{S^{n-1}}>0$. Since $\gamma$
is a diffeomorphism of $\mathbb R^m$ onto the upper half-sphere
$$
S^m_+=\{x=(x_1, x_2, \ldots, x_{m+1})\in\mathbb R^{m+1} :
|x|=1, ~x_{m+1}>0\}
$$
(where $m+1=n)$, we see that $|W|>0$.

   Thus we see that the special domain $G$ is generated by a pair
$(I, \varphi)$ where $\varphi\in C^{1, \omega}(I)$ is a function
with non-degenerate gradient and at the same time we have $1_G\in
A_p(\mathbb R^{m+1})$.

   By Lemma 1 we have
$$
\int_{\mathbb R} \frac{1}{|\lambda|^p}
\|e^{i\lambda\varphi}-1\|_{A_p(\mathbb R^m, I)}^p ~d\lambda<\infty,
$$
so
$$
\int_{\lambda\geq 1}\frac{1}{\lambda^p}
\|e^{i\lambda\varphi}-1\|_{A_p(\mathbb R^m, I)}^p d\lambda<\infty,
$$
and since $1\in A_p(\mathbb R^m, I), ~p>1,$ we see that
$$
\int_1^{\infty}\frac{1}{\lambda^p}
\|e^{i\lambda\varphi}\|_{A_p(\mathbb R^m, I)}^p d\lambda<\infty.
$$

  Hence, putting $V=I$ in estimate (4), we obtain
$$
\int_1^{\infty}\lambda^{m-p}
\bigg(\chi^{-1}\bigg(\frac{1}{\lambda}\bigg)\bigg)^{mp} d\lambda <\infty,
$$
that is (recall that $m=n-1$)
$$
\int_1^{\infty} \lambda^{n-1-p}
\bigg(\chi^{-1}\bigg(\frac{1}{\lambda}\bigg)\bigg)^{(n-1)p}
d\lambda <\infty.
$$

   The following lemma is of purely technical character.
It completes the proof of the theorem.

\quad

  \textbf{Lemma 2.} \emph{Let $n\geq 2, ~1<p<2$. The following
conditions are equivalent:}
$$
\int_1^{\infty} \lambda^{n-1-p}\bigg(\chi^{-1}
\bigg(\frac{1}{\lambda}\bigg)\bigg)^{(n-1)p}
d\lambda <\infty;
\leqno {1)}
$$
$$
\int_0^1 \frac{\delta^{n(p-1)-1}}{(\omega(\delta))^{n-p}}d\delta<\infty.
\leqno {2)}
$$

\quad

  Certainly, to complete the proof of the theorem it suffices to
verify that $1)\Rightarrow 2)$. The inverse implication for $n=2$
will be used below in $\S ~3$.

\quad

\emph{Proof of Lemma} 2. For $0<\varepsilon<1$ put
$$
I(\varepsilon)=\int_{1/\chi(1)}^{1/\chi(\varepsilon)}
\lambda^{n-1-p}\bigg(\chi^{-1}
\bigg(\frac{1}{\lambda}\bigg)\bigg)^{(n-1)p} d\lambda,
\qquad J(\varepsilon)=\int_\varepsilon^1
\frac{\delta^{n(p-1)-1}}{(\omega(\delta))^{n-p}}d\delta.
$$

   We have
$$
I(\varepsilon)=\frac{1}{n-p}\int_{1/\chi(1)}^{1/\chi(\varepsilon)}
\bigg(\chi^{-1}\bigg(\frac{1}{\lambda}\bigg)\bigg)^{(n-1)p} d\lambda^{n-p}.
$$
Changing the variable $\lambda=1/\chi(\delta)$ and integrating by
parts we obtain
$$
I(\varepsilon)=\frac{1}{n-p}
\bigg(\frac{\varepsilon^{n(p-1)}}
{(\omega(\varepsilon))^{n-p}}-
\frac{1}{(\omega(1))^{n-p}}\bigg)+\frac{(n-1)p}{n-p}J(\varepsilon).
\eqno(6)
$$

  Using this relation we see that
$$
I(\varepsilon)\geq \frac{-1}{(n-p)(\omega(1))^{n-p}}+
\frac{(n-1)p}{n-p}J(\varepsilon),
$$
so $1)\Rightarrow 2)$.

  Conversely, assume that condition 2) holds. Then, since
$$
\int_{\varepsilon/2}^\varepsilon
\frac{\delta^{n(p-1)-1}}{(\omega(\delta))^{n-p}}d\delta\geq
\frac{1}{(\omega(\varepsilon))^{n-p}}\int_{\varepsilon/2}^\varepsilon
\delta^{n(p-1)-1}d\delta\geq c_{n, p}\frac{\varepsilon^{n(p-1)}}
{(\omega(\varepsilon))^{n-p}},
$$
we have
$$
\frac{\varepsilon^{n(p-1)}}
{(\omega(\varepsilon))^{n-p}}\rightarrow 0,
\quad \varepsilon\rightarrow +0,
$$
and using (6) we obtain condition 1). The lemma and thus the
theorem are proved.

\quad

 \emph{Remark} 2. Theorem 2 (and Corollary 1, which it implies)
has local character. The theorem remains true if we assume that
only a part of the boundary of a domain $D$, i.e. the intersection
$B\cap\partial D$, where $B$ is a certain neighborhood in $\mathbb
R^n$, is $C^{1, \omega}$ -smooth and the normal map $\nu$ defined
on $B\cap\partial D$ is non-degenerate, that is $|\nu (B
\cap\partial D)|_{S^{n-1}}>0$. The condition that $D$ is bounded
can be replaced by the weaker condition $|D|<\infty$. (The
modification of the proof is obvious.)

  For $n=2$ the condition of non-degeneracy of the normal map on
$B\cap\partial D$ means that $B\cap\partial D$ is not a straight
line interval.

\quad

\begin{center}
\textbf{\S ~3. Domains in $\mathbb R^2$}
\end{center}

   In this section for each class $C^{1, \omega}$ (under certain
simple condition imposed on $\omega$) we shall construct a bounded
domain $D\subseteq\mathbb R^2$ with $C^{1, \omega}$ -smooth
boundary such that the characteristic function $1_D$ belongs to
$A_p$ for $p$ so close to $1$ as is allowed by Theorem 2. In
addition the domain $D$ has the property that its boundary does
not contain straight line intervals (thus, this domain is
essentially different from polygons).

  According to Theorem 2 if $D$ is a bounded domain in
$\mathbb R^2$ with $\partial D\in C^{1, \omega}$ and
$$
\int_0^1 \frac{\delta^{2p-3}}{\omega(\delta)^{2-p}} d\delta=\infty,
$$
then $1_D\notin A_p(\mathbb R^2)$. In particular this is the case
when $\partial D\in C^{1, \alpha}$ and $p\leq
1+\alpha/(2+\alpha)$. The following theorem shows that this result
is sharp.

\quad

\textbf{Теорема 3.} \emph{Suppose that
$\omega(2\delta)<2\omega(\delta)$ for all sufficiently small
$\delta>0$. There exists a bounded domain $D\subseteq\mathbb R^2$
with $\partial D\in C^{1, \omega}$ such that $1_D\in A_p(\mathbb
R^2)$ for all $p, ~1<p<2,$ satisfying
$$
\int_0^1 \frac{\delta^{2p-3}}{\omega(\delta)^{2-p}} d\delta<\infty.
\eqno(7)
$$
In addition the boundary of $D$ does not contain line intervals.}

\quad

This theorem immediately implies the following corollaries.

\quad

\textbf{Corollary 3.} \emph{For each $\alpha, ~0<\alpha<1,$ there
exists a bounded domain $D\subseteq\mathbb R^2$ such that its
boundary is $C^{1, \alpha}$ -smooth and $1_D\in A_p(\mathbb R^2)$
for all $p>1+\alpha/(2+\alpha)$. The boundary of $D$ does not
contain line intervals.}

\quad

\textbf{Corollary 4.} \emph{There exists a bounded domain
$D\subseteq\mathbb R^2$ such that its boundary is $C^1$ -smooth
and $1_D\in\bigcap_{p>1}A_p(\mathbb R^2)$. The boundary of $D$
does not contain line intervals.}

\quad

  Note also that from Theorems 2 and 3 it follows that the
existence of a domain $D\subseteq\mathbb R^2$ with $\partial D\in
C^{1, \omega}$ and $1_D\in\bigcap_{p>1} A_p(\mathbb R^2)$ is
equivalent to the condition that $\omega(\delta)$ tends to $0$
slower then any power, i.e., to the condition that
$\lim_{\delta\rightarrow+0} \omega(\delta)/\delta^\varepsilon
=\infty$ for all $\varepsilon>0$. Theorem 2 implies the necessity
of this condition. Theorem 3 implies its sufficiency.
\footnote{One should only observe that if this condition holds,
then one can find a nondecreasing continuous function
$\omega^\ast(\delta)$ on $[0, +\infty)$ that tends to $0$ slower
then any power and satisfies
$\omega^\ast(2\delta)<2\omega^\ast(\delta)$ for all $\delta>0$ and
$\omega^\ast(\delta)=O(\omega(\delta))$ as $\delta\rightarrow+0$.
For instance, one can put
$\omega^\ast(\delta)=\delta/(1+\delta)+\delta\inf_{0<x\leq \delta}
\omega(x)/x$. It is easy to verify that the second term is a
nondecreasing function, see [15, Ch. III, 3.2.5].}

  The author does not know whether similar results are true
for domains in $\mathbb R^n$ with $n\geq 3$.

\quad

\emph{Proof of Theorem} 3. Let $A_p(\mathbb T), ~1\leq
p\leq\infty,$ be the space of distributions $f$ on the circle
$\mathbb T=\mathbb R/2\pi\mathbb Z$ (where $\mathbb Z$ is the set
of integers) such that the sequence of Fourier coefficients
$\widehat{f}=\{\widehat{f}(k), ~k\in\mathbb Z\}$ belongs to $l^p$.
We put
$$
\|f\|_{A_p(\mathbb T)}=\|\widehat{f}\|_{l^p}=
\bigg(\sum_{k\in\mathbb Z}|\widehat{f}(k)|^p\bigg)^{1/p}.
$$
(For $1\leq p\leq 2$ each distribution in $A_p(\mathbb T)$ is a
function in $L^q(\mathbb T)\subseteq L^1(\mathbb T), ~1/p+1/q=1$.)

   For $p>1$ we put
$$
\Theta_p(y)=\bigg(\int_1^{y}
\bigg(\chi^{-1}\bigg(\frac{1}{\tau}\bigg)\bigg)^p d\tau \bigg)^{1/p},
\qquad y>1,
$$
where as above $\chi^{-1}$ is the function inverse to
$\chi(\delta)=\delta\omega(\delta)$.

   In Theorem 2 of the work [5], under assumption that
$\omega(2\delta)<2\omega(\delta)$ for all sufficiently small
$\delta>0$, we constructed a real function $\varphi$ on the circle
$\mathbb T=\mathbb R/ 2\pi\mathbb Z$ such that $\varphi\in C^{1,
\omega}(\mathbb T)$ (i.e., $\varphi$ is a $2\pi$ -periodic
function of class $C^{1, \omega}(\mathbb R)$) and for all $p,
~1<p<2,$ we have
$$
\|e^{i\lambda\varphi}\|_{A_p(\mathbb T)}\leq c_p \Theta_p(|\lambda|),
\qquad \lambda\in \mathbb R, \quad |\lambda|\geq 2.
\eqno(8)
$$
In addition the function $\varphi$ is nowhere linear, that is it
is not linear on any interval\footnote{Theorem 2 of the work [5]
contains similar result for $p=1$ as well.}.

  It is clear that from estimate (8) we have
$$
\|e^{i\lambda\varphi}-1\|_{A_p(\mathbb T)}\leq c_p \Theta_p(|\lambda|),
\qquad \lambda\in \mathbb R, \quad |\lambda|\geq 2.
\eqno(9)
$$

  It is also clear that for every continuously differentiable
function $f$ on $\mathbb T$ we have $f\in A_1(\mathbb T)$ and
$$
\|f\|_{A_1(\mathbb T)}\leq c \|f\|_{C^1(\mathbb T)},
$$
where
$$
\|f\|_{C^1(\mathbb T)}=
\max_{t\in \mathbb T}|f(t)|+\max_{t\in \mathbb T}|f'(t)|.
$$
So
$$
\|e^{i\lambda\varphi}-1\|_{A_p(\mathbb T)}\leq
\|e^{i\lambda\varphi}-1\|_{A_1(\mathbb T)}
$$
$$
\leq c \|e^{i\lambda\varphi}-1\|_{C^1(\mathbb T)}\leq
c_\varphi |\lambda|, \qquad \lambda\in\mathbb R.
\eqno(10)
$$

  Consider the following set $Q$ on the line $\mathbb R$:
$$
Q=\{t\in{(0, 2\pi) : \varphi'(t)>0}\}.
$$
Since $\varphi\neq \mathrm{const}$ and $\varphi(0)=\varphi(2\pi)$,
it is clear that $Q\neq\varnothing$ and $Q\neq (0, 2\pi)$.
Consider an interval $(a, b)$ which is a connected component of
the set $Q$. The derivative $\varphi'$ vanishes at least at one of
its endpoints. We can assume that at the right one, that is
$\varphi'(b)=0$, otherwise instead of $\varphi(t)$ and the
interval $(a, b)$ we consider the function $-\varphi(-t)$ and the
interval $(-b, -a)$. Choose now a point $c, ~a<c<b$. Replacing the
function $\varphi(t)$ by $\varphi(t)-\varphi(c)$, we can assume
that $\varphi(c)=0$. We put $I=(c, b)$.

  Thus we obtain a nowhere linear function
$\varphi\in C^{1, \omega}(\mathbb T)$ satisfying conditions (9),
(10) and an interval $I=(c, b)\subseteq [0, 2\pi]$ such that
$\varphi(c)=0$, the function $\varphi$ is strictly increasing on
$I$, and in addition $\varphi'(c)>0, ~\varphi'(b)=0$.

   Recall the well known relation [16, \S~44] between the spaces
$A_p(\mathbb T)$ and $A_p(\mathbb R)$ for $1<p\leq 2$. If $f$ is a
$2\pi$ -periodic function and $f^\ast$ is its restriction to $[0,
2\pi]$ extended by zero to $\mathbb R$, i.e. $f^\ast=f$ on $[0,
2\pi]$, $f^\ast=0$ on $\mathbb R\setminus[0, 2\pi]$, then $f\in
A_p(\mathbb T)$ if and only if $f^\ast\in A_p(\mathbb R)$. The
norms satisfy
$$
c_1(p)\|f^\ast\|_{A_p(\mathbb R)}\leq \|f\|_{A_p(\mathbb T)}
\leq c_2(p)\|f^\ast\|_{A_p(\mathbb R)}.
$$

   Thus, from estimates (9) and (10) we obtain that for
all $p, ~1<p<2,$
$$
\|e^{i\lambda\varphi}-1\|_{A_p(\mathbb R, I)}\leq
c_p \Theta_p(|\lambda|), \qquad
\lambda\in \mathbb R, \quad |\lambda|\geq 2, \eqno(11)
$$
and correspondingly
$$
\|e^{i\lambda\varphi}-1\|_{A_p(\mathbb R, I)}\leq
c|\lambda|, \qquad \lambda\in\mathbb R.
\eqno(12)
$$

   Consider the special domain $G\subseteq\mathbb R^2$
generated by the pair $(I, \varphi)$.

\quad

\textbf{Lemma 3.} \emph{For all $p, ~1<p<2,$ satisfying
\emph{(7)}, we have $1_G\in A_p(\mathbb R^2)$.}

\quad

\emph{Proof.} It is easy to verify that condition (7) implies
$$
\int_1^{\infty}\frac{1}{\lambda^p}(\Theta_p(\lambda))^p
d\lambda<\infty. \eqno(13)
$$
Indeed, for any $a>1$ integrating by parts we obtain
$$
\int_1^{a}\frac{1}{\lambda^p}(\Theta_p(\lambda))^p d\lambda=
\frac{1}{-p+1}\int_1^{a}(\Theta_p(\lambda))^p d\lambda^{-p+1}
$$
$$
=\frac{1}{-p+1}\bigg((\Theta_p(a))^p
a^{-p+1}-\int_1^{a} \lambda^{-p+1}
\bigg(\chi^{-1}\bigg(\frac{1}{\lambda}\bigg)\bigg)^p d\lambda \bigg)
$$
$$
\leq \frac{1}{p-1}\int_1^{a} \lambda^{-p+1}
\bigg(\chi^{-1}\bigg(\frac{1}{\lambda}\bigg)\bigg)^p d\lambda,
$$
and using Lemma 2 with $n=2$ we obtain (13).

  Therefore (see (11), (13)),
$$
\int_{|\lambda|\geq 2}\frac{1}{|\lambda|^p}
\|e^{i\lambda\varphi}-1\|_{A_p(\mathbb R, I)}^p d\lambda <\infty.
$$

  At the same time (see (12))
$$
\int_{|\lambda|<2}\frac{1}{|\lambda|^p}
\|e^{i\lambda\varphi}-1\|_{A_p(\mathbb R, I)}^p d\lambda <\infty.
$$
Thus,
$$
\int_{\mathbb R}\frac{1}{|\lambda|^p}
\|e^{i\lambda\varphi}-1\|_{A_p(\mathbb R, I)}^p d\lambda <\infty.
$$
It remains to use Lemma 1. The lemma is proved.

\quad

   Now we shall complete the proof of Theorem 3. The domain
$G\subseteq\mathbb R^2$, which we have constructed, is of the form
$$
G=\{(t,y): c<t<b, ~0<y<\varphi(t)\},
$$
where $\varphi\in C^{1, \omega}(\mathbb R)$ is a nowhere linear
function. Recall that according to our construction $\varphi(c)=0$
and the function $\varphi$ strictly increases on the interval $(c,
b)$. In addition $\varphi'(c)>0$ and $\varphi'(b)=0$. By Lemma 3
for all $p$ satisfying (7) we have $1_G\in A_p(\mathbb R^2)$.
Expanding (or contracting) the domain $G$ in the vertical
direction by an appropriate affine map, we can assume that
$\varphi'(c)=1$. Let $G^*$ be the domain symmetric to $G$ with
respect to the line $t=b$. Let $W=G\cup G^*\cup \xi$, where $\xi$
is the interval with the endpoints $(b, 0)$ and $(b, \varphi(b))$.
Take a square $\Pi\subseteq\mathbb R^2$ with side length $2(b-c)$.
We obtain the required domain $D$ by taking four rigid copies of
the domain $W$ (that is copies obtained by rotation and
translation) and gluing them to the sides of the square $\Pi$ on
its outer side. The theorem is proved.

\quad

\emph{Remark} 3. Theorem 3 (and its Corollaries 3, 4) allows the
following modification. The property that the boundary of the
domain $D$ does not contain line intervals can be replaced by the
property that $D$ is convex. The author does not know if it is
possible to get both this properties simultaneously. The indicated
modification follows since there exists a (real) nonconstant
function $\varphi\in C^{1, \omega}(\mathbb T)$ that satisfies (8)
and the following condition: the interval $(0, 2\pi)$ is a union
of three intervals such that the derivative $\varphi'$ is monotone
on each of them. Such a function is constructed by the author in
[5] (see the construction before the proof of Theorem 2 of [5]).

  It is not clear, even without any assumptions on smoothness
of the boundary, if there exists a strictly convex domain
$D\subseteq\mathbb R^n, ~n\geq 2,$ such that $1_D\in
\bigcap_{p>1}A_p(\mathbb R^n)$ (we call a domain strictly convex
if it is convex and its boundary does not contain line intervals).

\begin{center}
\textbf{References}
\end{center}

\flushleft
\begin{enumerate}

\item E. M. Stein, G. Weiss, \emph{Introduction to Fourier
    analysis on Euclidean spaces}, Princeton University Press,
    Princeton, New Jersy, 1971.

\item E. M. Stein, ``Problems in harmonic analysis related to
    curvature and oscillatory integrals'', \emph{Proceedings
    of the International Congress of Mathematicians, Berkely,
    CA, USA, 1986}, Amer. Math. Soc., 1987, pp. 196-221.

\item E. M. Stein, \emph{Harmonic analysis: Real-variable
    methods, orthogonality, and oscillatory integrals},
    Princeton University Press, Princeton, New Jersey, 1993.

\item C. S. Herz, ``Fourier transforms related to convex
    sets'', \emph{Ann. of Math.}, \textbf{75}:1(1962), 81-92.

\item V. V. Lebedev, ``Quantitative estimates in
    Beurling--Helson type theorems'', \emph{Sbornik:
    Mathematics}, \textbf{201}:12 (2010), 1811-1836.

\item V. V. Lebedev, ``Estimates in Beurling--Helson type
    theorems: Multidimensional case'', \emph{Mathematical
    Notes}, \textbf{90}:3 (2011), 373-384.

\item V. V. Lebedev, ``On the Fourier transform of the
    characteristic function of a domain in $R^n$'', 11-th
    Summer St. Petersburg Meeting in Math. Analysis; Euler
    Int. Math. Inst., St. Petersburg, Aug. 15-20, 2002,
    Abstracts, P. 24.

\item V. Lebedev, ``The domains with $C^1$ -smooth boundary
    and the Fourier transform of their characteristic
    functions'', 14-th Summer St. Petersburg Meeting in Math.
    Analysis; Euler Int. Math. Inst., St. Petersburg, June.
    6-11, 2005, Abstracts, P. 18.

\item V. Lebedev, ``On the Fourier transform of the indicator
    of a domain with $C^1$ -smooth boundary'', Int. Conference
    on Harmonic Analysis and Approximation III; Sept. 20-27,
    2005, Tsahkadzor, Armenia, Abstracts, P. 50-51, Yerevan,
    2005.

\item H. Triebel, ``Function spaces in Lipschitz domains and
    on Lipschitz manifolds. Characteristic functions as
    pointwise multipliers'', \emph{Rev. Mat. Complutense},
    \textbf{15} (2002), 475-524.

\item E. M. Stein, \emph{Singular integrals and
    differentiability properties of functions}, Princeton
    University Press, Princeton, New Jersy, 1970.

\item M. N. Kolountzakis, T. Wolff, ``On the Steinhaus tiling
    problem'', \emph{Mathematika}, \textbf{46} (1999),
    253-280.

\item P. Mattila, \emph{Geometry of sets and measures in
    Euclidean spaces}, Cambrige Univ. Press, 1995.

\item L. H\"ormander, ``Estimates for translation invariant
    operators in $L^p$ spaces'', Acta Math. \textbf{104}
    (1960), 93-140.

\item A. F. Timan, \emph{Theory of approximation of functions
    of a real variable}, Pergamon Press, Oxford--London--New
    York--Paris, 1963.

\item  M. Plancherel, G. P\'olya, ``Fonctions enti\`eres et
    int\'egrales de Fourier multiples. II'', \emph{Comment.
    Math. Helv.}, \textbf{10}:2 (1937), 110-163.

\end{enumerate}

\quad

\qquad\textsc{V. V. Lebedev}\\
\qquad Dept. of Mathematical Analysis\\
\qquad Moscow State Institute of Electronics\\
\qquad and Mathematics(Technical University)\\
\qquad E-mail address: \emph {lebedevhome@gmail.com}

\end{document}